\documentclass[12pt]{article}
\usepackage{cancel}
\usepackage[usenames]{color}
\usepackage{amssymb}
\usepackage{amscd}
\usepackage[colorlinks=true,linkcolor=webgreen,
filecolor=webbrown,
citecolor=webgreen]{hyperref}

\definecolor{webgreen}{rgb}{0,.5,0}
\definecolor{webbrown}{rgb}{.6,0,0}

\usepackage{color}
\usepackage{fullpage}
\usepackage{float}

\usepackage{amsmath}
\usepackage{amsthm}
\usepackage{amsfonts}

\newtheorem{theorem}{Theorem}
\newtheorem{lemma}{Lemma}
\newtheorem{problem}{Open Problem}
\newtheorem{conjecture}{Conjecture}
\theoremstyle{definition}\newtheorem{remark}{Remark}
\title{Words with factor complexity $2n+1$ and minimal critical exponent}
\author{James D. Currie\\
Department of Mathematics and Statistics\\ 
University of Winnipeg\\
Winnipeg, Manitoba R3B 2E9, Canada}

\begin{document}
\maketitle
\begin{abstract} Let the word ${\mathbf G}$ be the fixed point of the morphism $\gamma$ sending $0$ to $01$, $1$ to $2$, and $2$ to $02$.
In 2019, Shallit and Shur showed that ${\mathbf G}$ has factor complexity $2n+1$. They also showed that ${\mathbf G}$ has critical exponent $\mu=2+\frac{1}{\lambda^2-1}= 2.4808726\cdots$, where $\lambda=1.7548777$ is the real zero of $x^3-2x+x-1=0$. They conjectured that this was the least possible critical exponent among the words with factor complexity $2n+1$.
We confirm their conjecture. We anticipate that our method, including an intricate case analysis by computer, will have wider application.\vspace{.1in}


\end{abstract}

\section{Introduction}

The critical exponent of a word is the supremum of the exponents of its factors. (Definitions of basic terms in combinatorics on words, such as factors, exponents, etc., follow in the next section.) In recent years a large number of papers have been written answering questions of the following form:\\

 Given a class ${\cal L}$ of infinite words, what is the minimum critical exponent of a member of ${\cal L}$?\\

The various steps in the solution of Dejean's conjecture by several authors \cite{dejean,pansiot,ollagnier,morteza,carpi,rampersad1,rao} answered this question for various classes of the form ${\cal L}=\Sigma_k^\omega$ where $\Sigma_k=\{0,1,2,\ldots,k-1\}$.
Since the resolution of Dejean's conjecture in 2011 \cite{rampersad1,rao} researchers have explored versions of this question 
where ${\cal L}$ is variously taken  to be, e.g., the class of
Sturmian words \cite{luca},
rich words over certain alphabets \cite{mol1,mol2},
balanced words \cite{rampersad2,dvorakova1},
  or 
complementary symmetric Rote words \cite{dvorakova2}.

For example, the minimum critical exponent for a Sturmian word is $$2+\varphi=3.618033988749\cdots$$ where $\varphi$ is the golden ratio. This is realized by the Fibonacci word, the fixed point of the morphism $\Phi$ sending $0$ to $01$ and $1$ to $0$. One starting point for the study of Sturmian words and their generalizations is the relevant chapter in Allouche and
Shallit's book \cite{allouche}.

Words of factor complexity $2n+1$ are considered in the paper of Cassaigne, Labb\'{e}, and Leroy \cite{cassaigne}. Such words are connected to the well-studied Arnoux-Rauzy words: the Arnoux-Rauzy words over a $k$-letter alphabet have factor complexity $(k-1)n+1$. Note that the Arnoux-Rauzy words over a binary alphabet are just the Sturmian words.  In the present paper, we find the minimal critical exponent for the class of words with
factor complexity $2n+1$. This confirms the 2019 conjecture of 
Shallit and Shur \cite{shallit2019}. Taken together with 
Carpi and De Luca's result \cite{luca} for Sturmian words, this gives the beginnings of a Dejean-type theorem for Arnoux-Rauzy words: we know the minimal critical exponent in the cases $k=2$ and $k=3$.

To find the minimal critical exponent for binary rich words, it was necessary to prove the following structure theorem:

\begin{theorem}\cite{mol1}
Let ${\mathbf w}\in\Sigma_2^\omega$ be a $14/5$-power-free rich word.  For every $n\geq 1$, a suffix of ${\mathbf w}$ has the form $f(\phi^n(w_n))$ or $f(g(\phi^n(w_n)))$ for some word $w_n\in\Sigma_3^\omega$ where
\begin{align*}f(0)=&0&g(0)=&011&\phi(0)=&01\\
f(1)=&01&g(1)=&0121&\phi(1)=&02\\
f(2)=&011&g(2)=&012121&\phi(2)=&022.
\end{align*}
\end{theorem}

Paraphrasing, if we let ${\mathbf H}$ be the fixed point of $\phi$, then the infinite $14/5$-power-free rich binary words look like one of $h_1({\mathbf H})$ and $h_2({\mathbf H})$, where $h_1=f$, and $h_2=f\circ g$.

In the present paper, we prove the analogous structure theorem for 
the infinite words with factor complexity $2n+1$, and not containing $r$-powers for $r\ge 5/2$. However, instead of the two possible morphisms $h_1$ and $h_2$ in the above-mentioned theorem, there are several possible morphisms. (See Figures~\ref{morphisms} and \ref{dual morphisms} and following.)
The case analysis is therefore complex. (Interestingly, Shallit and Shur \cite{shallit2019} prove that the minimum critical exponent for the infinite words with factor complexity $2n$ is exactly $5/2$.) 

For several decades, the use of computer backtracking has been a standard tool in combinatorics on words. More recently, however, this  area has seen a number of trends where computer-assisted proofs come closer to being `computer proofs'. One such trend has been seen in the continued success of the Walnut computer  package \cite{mousavi,shallit2022}, which proves properties of automatic sequences. Over 100 papers using Walnut are listed by Shallit \cite{walnut}. Of relevance to the current paper, in many situations Walnut can check the critical exponent of the fixed point of a morphism, as in a recent example by Baranwal et al. \cite{baranwal}. 

Another instance of the intensive use of computation in proofs has been seen in papers using the `template method' or its variations \cite{andrade,dekking,template,rosenfeld1,rosenfeld2,cassaigne}, where one
recursively generates ancestors of patterns of interest that could potentially appear in
the fixed point of some morphism.

The present paper presents what is, to our knowledge, a new 
 sort of application of computation to combinatorics on words. In  analyzing the subject of this paper, a great proliferation of subcases arose, branching to a depth of, for example, subcase A.1.2.2.2 (5 levels of subcases). Some cases involved branching into as many as 7 subcases. Each subcase involved backtrack searches and might have one of several resolutions, possibly including the spawning of new subcases.

Having started with a case analysis by hand, at a certain point we realized that the spawning and resolution of subcases would be best handled by computer. In the end, we wrote a computer program which in turn wrote the desired proof.
\section{Word preliminaries}
An {\em alphabet} is a finite set with elements known as {\em letters}. In what follows, we use the alphabet $\Sigma_3=\{0,1,2\}$.
The set of finite words over $\Sigma_3$ is denoted by 
$\Sigma_3^*$, and can be regarded as the free semigroup generated by the letters $0$, $1$, and $2$. The length of a word $u$ is the number of occurrences of letters in it and is denoted by $|u|$. Thus, for example,
$|01202|=5$. For each $a\in\Sigma_3$, we denote the number of occurrences of letter $a$ in $u$ by $|u|_a$. Thus, for example,
$|01202|_0=2,|01202|_1=1,|01202|_2=2$.

If $u$, $v\in\Sigma_3^*$, the {\em concatenation} $uv$ of $u$ and $v$ consists of the letters of $u$ followed by the letters of $v$. Thus if $u=012$ and $v=02$, we have $uv=01202$. We say that $u$ is a {\em prefix} of $uv$ and $v$ is a suffix of $uv$. We also say that $uv$ is a {\em right extension} of $u$ {\em by} $v$. If $u$, $v$, $w\in\Sigma_3^*$, we say that $v$ is a {\em factor} of 
$uvw$. If word $y$ is not a factor of word $z$ we say that $z$ {\em omits} the factor $y$. Thus $01$, $120$, and $02$ are respectively a prefix, factor, and suffix of $01202$, and $01202$ is a right extension of $01$ by $202$. Word $01202$ omits the factor $10$. 

Suppose $u=u_0u_1u_2\cdots u_{n-1}$ where each $u_i\in\Sigma_3$. We say that positive integer $p$ is a {\em period} of $u$ if
$u_i=u_{i+p}$ whenever $0\le i\le i+p\le n-1$. In this case we say that $u$ is an {\em $r$ power}, where $r=n/p$. Thus 3 is a period of the word $01201201$, which is an $8/3$ power. 
Call a word {\em low} if it  contains no  $r$ power for  $r\ge 5/2$.

A {\em morphism} on $\Sigma_3$ is a semigroup homomorphism of $\Sigma_3^*$. A morphism is determined by its images on letters. Of particular interest is the morphism $\gamma:\Sigma_3^*\rightarrow\Sigma_3^*$ generated by
\begin{align*}
\gamma(0)&=01\\
\gamma(1)&=2\\
\gamma(2)&=02
\end{align*}
For conciseness we record a morphism $h$ as $h=[h(0),h(1),h(2)]$.
Thus $\gamma=[01,2,02]$.
If $|h(a)|>0$ for each $a\in\Sigma_3$, we say that $h$ is {\em non-erasing}. We use exponentiation to denote iteration of a morphism. Thus $h^2(u)=h(h(u))$, $h^3(u)=h(h(h(u)))$, etc.

Suppose $u=u_0u_1u_2\cdots u_{n-1}$ where each $u_i\in\Sigma_3$. The {\em reverse} of $u$ is the word
$u^R=u_{n-1}\cdots u_2u_1u_0$. The reverse of the morphism $h=[h(0),h(1),h(2)]$ is defined to be $h^R=[(h(0))^R,(h(1))^R,(h(2))^R]$. Thus, for example, 
$\gamma^R=[10,2,20]$. We see that for a word $u$ and the morphism $h$ we have $(h(u))^R=h^R(u^R).$

As well as the finite words, we consider the right-infinite words over $\Sigma_3$. Formally, such a word is a sequence 
${\mathbf u}=\{u_n\}_{n=0}^\infty$ where each $u_i$ is a single letter.
We use bold-face letters to distinguish right-infinite words from finite words. For each non-negative integer $n$, the length-$n$ prefix of ${\mathbf u}$ is the finite word $u_0u_1u_2\cdots u_{n-1}$. We say that $v\in\Sigma_3^*$ is a factor of ${\mathbf u}$ if it is a factor of some finite prefix of ${\mathbf u}$. Call a right-infinite word {\em low} if each of its finite prefixes is low.

We extend the notaion of suffix to right-infinite words; a suffix
of ${\mathbf u}$ is a right-infinite word 
${\mathbf u}'=\{u_{n+i}\}_{n=0}^\infty$ for some non-negative
integer $i$. We write ${\mathbf u}=u{\mathbf u}'$ where $u$ is the prefix of ${\mathbf u}$ of length $i$.

Let ${\mathbf w}\in\Sigma_3^\omega$. Fix a finite factor $b$ of ${\mathbf w}$. The word $bu$ is a {\em $b$-block} of ${\mathbf w}$ if 
\begin{itemize}
\item The word $bub$ is a factor of ${\mathbf w}$;
\item Whenever $p$ is a prefix of $bub$ such that $b$ is a suffix of 
$p$, then $p=b$ or $p=bub$.
 \end{itemize}
If $b$ appears infinitely often as a factor of ${\mathbf w}$ then a suffix of ${\mathbf w}$ can be concatenated from the $b$-blocks of ${\mathbf w}$.

The set of all right-infinite words over $\Sigma_3$ is denoted by $\Sigma_3^\omega$. Right-infinite words arise in various ways: If $u=u_0u_1\cdots u_{p-1}$ where $p$ is positive and the $u_i$ are letters, then $u^\infty=\{u_i\}_{i=0}^\infty$ where $u_i=u_{i \bmod p}$. Thus, for example,
$$(012)^\omega=012012012012\cdots.$$
Another way a right-infinite word arises is via certain morphisms. Let $h$ be a morphism such that $a$ is a prefix of $h(a)$ for some
$a\in\Sigma_3$. By induction, $h^{n-1}(a)$ is a prefix of $h^n(a)$ for each $n$. If $|h^n(a)|$ increases without bound, we can define 
$h^\omega(a)$ to be the unique right-infinite word having  
$h^n(a)$ as a prefix for each $n$.  For a right-infinite word ${\mathbf u}$ over $\Sigma_3$, we define $h({\mathbf u})$ to be the right-infinite word having $h(u)$ as a prefix for each prefix $u$ of  ${\mathbf u}$; thus  $h^\omega(a)$ is a fixed point of $h$.
In this paper we are particularly interested in the fixed point
${\mathbf G}=\gamma^\omega(0)$.

A third (non-constructive) way to obtain right-infinite words is via K\"onig's lemma \cite{kleene}. Here is one formulation of K\"onig's lemma in  the combinatorics on words context: A subset $L\subseteq\Sigma_3^*$ is said to be a {\em factorial language} if whenever $\ell$ is in $L$, every factor of 
$\ell$ is also in $L$.
\begin{lemma}[K\"onig's lemma] Let $L\subseteq\Sigma_3^*$ be an infinite factorial language. Then there is a right-infinite word
${\mathbf u}\in\Sigma_3^\infty$ such that every prefix of ${\mathbf u}$ is in $L$.
\end{lemma}

If ${\mathbf u}\in\Sigma_3^\omega$, the {\em factor complexity} of ${\mathbf u}$ is the function
counting the length $n$ factors of ${\mathbf u}$. Thus, for example, if ${\mathbf u}=(012)^\omega$, then the factor complexity of ${\mathbf u}$ is $C(n)$ where
\begin{align*}
C(0)&=1\\
C(n)&=3, n\ge 1.
\end{align*}
Shallit and Shur \cite{shallit2019} have shown that ${\mathbf G}$ has factor complexity $2n+1$.

We say that ${\mathbf w}\in\Sigma_3^\omega$ is {\em eligible} if it is low  and has factor complexity $2n+1$. 

The {\em critical exponent} of a finite or right-infinite word is the supremum of those $r$ such that a factor of the word is an $r$ power. Shallit and Shur \cite{shallit2019} have shown that ${\mathbf G}$ has critical exponent $\mu=2+\frac{1}{\lambda^2-1}= 2.4808726\cdots$, where $\lambda=1.7548777$ is the real zero of $x^3-2x+x-1=0$. They conjecture that this is the least possible critical exponent among words with factor complexity $2n+1$.

We say that the right-infinite word ${\mathbf u}\in\Sigma_3^\omega$ has {\em letter frequencies} 
$\rho_0$, $\rho_1$, $\rho_2$, if for all $\epsilon_1>0$ there is an $N_1$ such that if $u$ is a factor ${\mathbf u}$ of length at least $N_1$, then \begin{equation}\rho_a-\epsilon_1 <\frac{|u|_a}{|u|}< 
\rho_a+\epsilon_1,\text{ for }a\in\Sigma_3.\end{equation}

Using standard methods, it can be shown that  ${\mathbf G}$ has letter frequencies.
(For values of the $\rho_a$ see the entry for
sequence \href{https://oeis.org/A287104}{A287104} in the On-Line Encyclopedia of Integer Sequences; the word ${\mathbf G}$ is a recoding of this sequence.)

\section{Main Theorem}

\begin{theorem}[Main Theorem]\label{main theorem} Let ${\mathbf w}\in\Sigma_3^\omega$ have  
factor complexity $2n+1$. Then the critical exponent of ${\mathbf w}$ is at least $\mu=2+\frac{1}{\lambda^2-1}= 2.4808726\cdots$, where $\lambda=1.7548777$ is the real zero of $x^3-2x+x-1=0$. 
\end{theorem}

A key ingredient in the proof of this main theorem is a  structure theorem. 
\begin{theorem}[Structure Theorem]\label{structure theorem} Let ${\mathbf w}\in\Sigma_3^\omega$ be eligible. There is a non-erasing morphism $h$ such that either
\begin{itemize}
\item 
For every factor $g$ of ${\mathbf G}$, the word $h(g)$ is a factor of ${\mathbf w}$.
\item 
For every factor $g$ of ${\mathbf G}$, the word $h(g^R)$ is a factor of ${\mathbf w}$.
\end{itemize}
\end{theorem}

The main theorem follows from this lemma:

\begin{lemma}\label{nu}
Let ${\mathbf u}\in\Sigma_3^\omega$ have irrational critical exponent $\nu$, with $2<\nu<3.$ 
Suppose that  ${\mathbf u}$ has letter frequencies $\rho_0$, $\rho_1$, $\rho_2$.
Let ${\mathbf w}\in\Sigma_3^\omega$ and let $h:\Sigma_3^*\rightarrow \Sigma_3^*$ be a non-erasing morphism. Suppose that for each factor $g$ of ${\mathbf u}$, either $h(g)$ or $h(g^R)$ is a factor of ${\mathbf w}$. Then the critical exponent of ${\mathbf w}$ is at least $\nu$.
\end{lemma}

\begin{remark} The use of the alphabet $\Sigma_3$ and the restriction that $2<\nu<3$ are not essential to this lemma, but simplify the notation in its proof.
\end{remark}

\begin{proof}[Proof of main theorem] Let ${\mathbf w}\in\Sigma_3^\omega$ have  
factor complexity $2n+1$. Recall that $\mu=2+\frac{1}{\lambda^2-1}= 2.4808726\cdots$, where $\lambda=1.7548777$ is the real zero of $x^3-2x+x-1=0$.
Thus $\mu<5/2$. If ${\mathbf w}$ is not low, then the critical exponent of ${\mathbf w}$ is at least $5/2>\mu$, and we are done.

Suppose then that ${\mathbf w}$ is low. By the structure theorem, there is a non-erasing morphism $h$ such that for every factor $g$ of ${\mathbf G}$, either $h(g)$ is a factor of ${\mathbf w}$ or $h(g^R)$ is a factor of ${\mathbf w}$. (This is actually weaker than the conclusion of the structure theorem).
The word  ${\mathbf G}$ has letter frequencies.
It is easy to show that $\mu$ is irrational.
By  Lemma~\ref{nu}, then,  the critical exponent of ${\mathbf w}$ is at least $\mu$.
\end{proof}
\section{Eligible words and their factors}
We defer the proof of the following lemma.
\begin{lemma}\label{0,1,2} Suppose ${\mathbf w}\in\Sigma_3^\omega$ is eligible. Then every suffix of ${\mathbf w}$ contains letters 0, 1, and 2.
\end{lemma}
Suppose that ${\mathbf w}$ is eligible. Consider the directed graph $D$ with vertex set $\Sigma_3$, and with
a directed edge $xy$ for $x,y\in\Sigma_3$ exactly when $xy$ is a factor of ${\mathbf w}$. Then there is a directed path on $D$ labeled with the letters of ${\mathbf w}$. Since ${\mathbf w}$ has exactly 5 length-2 factors, the graph $D$ has exactly 5 edges.

By the previous lemma, the graph $D$ must be strongly connected; that is, there is a directed path between any two vertices of $D$.  Therefore, $D$ must have at least 3 non-loop edges, and if it has exactly three non-loop edges, they form a directed cycle.

It is convenient to keep track of words which are not factors of ${\mathbf w}$: Let  $F=\{u\in\Sigma_3^*:u\text{ is not a factor of }{\mathbf w}\}.$ To prove the structure theorem, it suffices to consider two cases involving $F$.

\begin{itemize}
\item If $D$ has exactly three non-loop edges (forming a directed cycle) then, permuting $\Sigma_3$ if necessary, we may assume that
$$\{01,12,20,00\}\subseteq F.$$
\item If $D$ has more than three non-loop edges, permuting the alphabet if necessary,we may assume that
$$\{11,22\}\subseteq F.$$
\end{itemize}

The following technical lemmas give sufficient conditions for concluding that the structure theorem holds:

\begin{lemma}[Morphism Lemma]\label{morphism lemma}
Suppose that ${\mathbf w}\in\Sigma_3^\omega$ is low.  Suppose that some suffix of ${\mathbf w}$ has the form $h({\mathbf u})$, some ${\mathbf u}\in\Sigma_3^\omega$ where $h=[a,b,c]$, some $a,b,c\in\Sigma_3^+$ such that
\begin{enumerate}
\item The word $b$ is a  prefix of $a$, which is a  prefix of $c$.
\item We have $|b|\ge |c|/2$. 
\item The words $b$ and $c$ have a common suffix $s$ such that $|s|\ge |b|/2$.
\end{enumerate}
If $g$ is any factor of ${\mathbf G}$, then $h(g)$ is a factor of ${\mathbf w}$.
\end{lemma}
\begin{lemma}[Dual Morphism Lemma]\label{dual morphism lemma}
Suppose that ${\mathbf z}\in\Sigma_3^\omega$ is low.  Suppose that a suffix of ${\mathbf z}$ has the form $h({\mathbf u})$ for some ${\mathbf u}\in\Sigma_3^\omega$, where $h=[a,b,c]$ for some $a,b,c\in\Sigma_3^+$ such that
\begin{enumerate}
\item The word $b$ is a  suffix of $a$, which is a  suffix of $c$.
\item We have $|b|\ge |c|/2$. 
\item The words $b$ and $c$ have a common prefix $p$ such that $|p|\ge |b|/2$.
\end{enumerate}
If $g$ is any factor of ${\mathbf G}$, then $h(g^R)$ is a factor of ${\mathbf z}$.
\end{lemma}

\begin{lemma}[First Parsing Lemma]\label{first parsing lemma}
Let $S\subseteq\Sigma_3^*$ and let $b\in\Sigma_3^*$. Let $d$ be a non-negative integer. Suppose that every low word of length $d+1$ contains a factor from $S\cup\{b\}$. Then if ${\mathbf w}$ is an eligible right-infinite word which omits the words of $S$, a suffix of ${\mathbf w}$ can be concatenated from $b$-blocks of the form $bu$, where $|u|\le d$. \end{lemma}

\begin{lemma}[Second Parsing Lemma]\label{second parsing lemma} Let $S,B\subseteq\Sigma_3^*$. Let ${\mathbf w}$ be an eligible right-infinite word which omits the words of $S$, such that ${\mathbf w}$ is concatenated from words of $B$. Let $s\in B$ have the property that for all $r,t\in B$, either $rst$ is not low, or $rst$ contains a factor from $S$. Then a suffix of ${\mathbf w}$  is concatenated from words of $B-\{s\}$.
\end{lemma}

\section{Proof of the structure theorem}
We prove the structure theorem via an intricate case analysis. Each case is determined by a set $S\subseteq F$, where $F$ is the set of ternary words which are not factors of ${\mathbf w}$. As remarked in the previous section, it suffices to resolve the subcases $S=\{01,12,20,00\}$ and $S=\{11,22\}$.

In a given case, labelled by a set $S$, we proceed in one of three ways. We either
\begin{enumerate}
\item Show that the conditions of the morphism lemma or the dual morphism lemma hold, so that the case is resolved;
\item Show that the case cannot arise, so that the case is resolved;
\item Partition the case into subcases.
\end{enumerate}

The proof is thus carried out recursively, in a depth-first way. In situation 1 or 2, a case is resolved. In situation 3, new subcases are added. If at some point all subcases have been resolved (and Hercules beats the hydra) the structure theorem has been proved.

Let a set $S$ of forbidden factors be given. Say that a word $w$ is {\em good} if it is low and has no factor in $S$.
The concrete recipe for carrying out 1, 2, and 3 above is as follows.
\begin{enumerate}
\item[A.] If $T\subseteq S$, and case $T$ has been previously resolved, then $S$ is resolved. 
\item[B.] Perform a backtrack search, looking for a good word $w$ of length 250. Word $w$ is a hypothetical prefix of ${\mathbf w}$.  
If no such word $w$ exists, the case cannot arise.
\item[C.] If we find $w$ of length 250, count the factors of $w$ of each length $n$ with $1\le n\le 20$, until an $n$ is found such that $w$ has more than $2n+1$ factors of length $n$. 
\begin{enumerate}
\item[a.] If we find $n$ such that $w$ has more than $2n+1$ factors of length $n$, then we classify each such factor $s$ as {\em needed} or {\em unneeded}. A factor $s$ is needed if a backtrack shows that no low word of length 250 omits all the factors of $S\cup\{s\}$; it is unneeded otherwise.
\begin{enumerate}
\item[i.] If there are more than $2n+1$ needed factors then no low word of length 250 can be a prefix of ${\mathbf w}$, since the condition on the complexity of ${\mathbf w}$ is violated. The case where ${\mathbf w}$ omits factors in $S$ cannot arise.
\item[ii.] Otherwise, for each unneeded factor $s$ we form a subcase, replacing $S$ by $S\cup\{s\}$.
\end{enumerate}  
\item[b.]  If no such $n$ is found, we look for a morphism satisfying the conditions of the morphism lemma or the dual morphism lemma. To do this
\begin{enumerate} 
\item Consider the good words $b$ with $|b|\le 3$.
\item Suppose that for some positive integer $d\le 200$, every good word of length $d+1$ has $b$ as a factor. This can be checked by a finite search.
\item Then consider the set $B=\{bu\in\Sigma_3^*:bub$ is good, $|u|\le d$, and whenever $p$ is a prefix of $bub$ such that $b$ is a suffix of $p$, then $p=b$ or $p=bub\}$. By the first parsing lemma, ${\mathbf w}$ is concatenated from words of $B$.
\item We form a morphism using the blocks in B if either:
\begin{itemize}
\item Set $B$ contains exactly 3 $b$-blocks;
\item Set $B$ contains exactly 4 $b$-blocks, including a $b$-block $s$ such that for all $r,t\in B$, word $rst$ is not good. In this case we replace $B$ by $B-\{s\}$. By the second parsing lemma, a final segment of ${\mathbf w}$ is concatenated from the $b$-blocks of $B$.
\end{itemize}
\end{enumerate}

We form a morphism from the three $b$-blocks of $B$, and consider this morphism and its conjugates.
If one of the considered morphisms satisfies the conditions of the morphism lemma or the dual morphism lemma the case is resolved.
\end{enumerate}
  \end{enumerate}

\begin{remark} 
The `hard-wired' constants 250 (for backtrack search), 3 (for $|b|$), 20 (depth of the complexity check) are all {\em ad hoc}. They were originally taken to be larger, but then tuned.
\end{remark}
\begin{proof}[Proof of the structure theorem]
A python implementation of the above scheme in SageMath 9.3 (on a 
Microsoft Surface with Intel(R) Core(TM) i5-1035G4 CPU at 1.10GHz)
resolved the case $S=\{11,22\}$ in about 36 seconds, and the case $S=\{01,12,20,00\}$ in about 9 seconds. The code is available at
\url{https://github.com/drjamescurrie-ai/RACE-3/tree/main}.
\end{proof}
\begin{remark} In fact, the implementation mentioned above kept track of sets $S$ and the reasoning at each of its steps, and produced `human readable' proofs. We give the reasoning that was output for  the case $S=\{01,12,20,00\}$ in an appendix. Further, the implementation kept track of all the morphisms $h$ used. It may be the case that for some of these morphisms, the word $h({\mathbf G})$ is eligible and has critical exponent $\mu$.
\end{remark}

\section{Proofs of Lemmas}
\begin{proof}[Proof of Lemma~\ref{nu}]
The critical exponent of ${\mathbf u}$ is irrational, and therefore is not realized by any of its finite factors. Thus for all
$\epsilon_2>0$ and any positive integer $N$, the word ${\mathbf u}$ has a factor $uuu'$ such that $u'$ is a prefix of $u$ with $|u'|>N$ and \begin{equation}\label{epsilon_2}\frac{|uuu'|}{|u|}>\nu-\epsilon_2.\end{equation} 

Since $h$ is non-erasing, we have
$$\sum_{a=0}^2|h(a)|\rho_a\ge\sum_{a=0}^2\rho_a=1$$ and is non-zero. Then
$$\lim_{\epsilon\rightarrow 0^+} \frac{\sum_{a=0}^2|h(a)|(\rho_a-\epsilon)}{\sum_{a=0}^2|h(a)|(\rho_a+\epsilon)}=1$$

Given $\epsilon>0$, choose $\epsilon_1>0$ such that
$$\frac{\sum_{a=0}^2|h(a)|(\rho_a-\epsilon_1)}{\sum_{a=0}^2|h(a)|(\rho_a+\epsilon_1)}>1-\epsilon.$$

Choose $N_1$ such that when $|u'|>N_1$, then 
$$\rho_a-\epsilon_1 <\frac{|u|_a}{|u|}< 
\rho_a+\epsilon_1,\text{ for }a\in\Sigma_3.$$
 Given $\epsilon_2>0$, choose a factor $g=uuu'$ of ${\mathbf u}$ with
$|u'|>N_1$ so that (\ref{epsilon_2}) holds. 
(If $h(g^R)$ is a factor of ${\mathbf w}$, we replace $g$ by $g^R$, redefining $u$ and $u'$ accordingly.)

 The word
$h({\mathbf u})$ contains the factor
$h(uuu')=h(u)h(u)h(u')$, and $h(u')$ is a prefix of $h(u)$. 

Also
\begin{align*}
\frac{|h(uuu')|}{|h(u)|}
&=\frac{\sum_{a=0}^2|h(a)||uuu'|_a}{\sum_{a=0}^2|h(a)||u|_a}\\
&\ge \frac{\sum_{a=0}^2|h(a)|(\rho_a-\epsilon)|uuu'|}{\sum_{a=0}^2|h(a)|(\rho_a+\epsilon)|u|}\\
&= \frac{|uuu'|}{|u|}\frac{\sum_{a=0}^2|h(a)|(\rho_a-\epsilon)}{\sum_{a=0}^2|h(a)|(\rho_a+\epsilon)}\\
&>(\nu-\epsilon_2)(1-\epsilon).\end{align*} 
We see that ${\mathbf w}$ contains factors with exponent arbitrarily close to $\nu$. The result follows.
\end{proof} 
\begin{proof}[Proof of Lemma~\ref{0,1,2}]
Let ${\mathbf s}$ be a suffix of ${\mathbf w}$. It suffices to show that ${\mathbf s}$ contains letter 0. 

Suppose ${\mathbf s}$ contains only letters 1 and
2. It must contain factors 12 and 21; otherwise ${\mathbf s}$ would end with  $1^\omega$ or $2^\omega$ and therefore contain $5/2$ powers.

Suffix ${\mathbf s}$ must also contain factor 11 (and, symmetrically, factor 22): If not, assume ${\mathbf s}$ starts with 1, replacing ${\mathbf s}$ with one of its suffixes if necessary. Then ${\mathbf s}$ is concatenated from the $1$-blocks of ${\mathbf w}$. Since 222
is a $3$ power,
the only possible $1$-blocks of ${\mathbf w}$ are  12 and 122. The $1$-block 12 can only be used once in this concatenation, or else it appears in the context $2\underline{12}12$, which is a $5/2$ power.  Then a suffix is $(122)^\omega$, which contains  $5/2$ powers. This is a contradiction. Thus 11 must indeed occur.

So far we have shown that ${\mathbf s}$ contains factors 11, 12, 21, and 22. Since ${\mathbf w}$ contains $2(2)+1=5$ factors of length 2, there is exactly one length-2 factor of ${\mathbf w}$ containing a $0$. It follows that the only 0 in ${\mathbf w}$ is its first letter. We thus see that ${\mathbf w}$ contains no factor of the form $a0$ where $a\in\Sigma_3$.

Consider the following set of 10 length-4 words:
$$S=\{ 1121, 1122, 1211, 1212, 1221, 2112, 2121, 2122, 2211, 2212\}.$$

Let $s$ and $t$ be distinct words of $S$. Suppose $u$ is a low word over $\Sigma_3$ such that the set of factors of $u$ does not include any of 00, 10, 20, $s$, or $t$. 
For each of the $\binom{10}{2}$ possibilities for $s$ and $t$, a backtrack search shows that the length of $u$ is at most 39.  
It then follows that ${\mathbf w}$ contains 9 out of 10 
words of $S$ as factors. In addition,the word  ${\mathbf w}$ contains a length-4 factor starting with 0. However, the word ${\mathbf w}$ is eligible, and cannot contain $10$  length-4 factors. This is a contradiction.
\end{proof}

We use two preliminary lemmas in the proof of the morphism lemma (Lemma~\ref{morphism lemma}).
\begin{lemma}\label{00,11,21,22}
Suppose that ${\mathbf w}\in\Sigma_3^\omega$ is low. Let $a,b,c\in\Sigma_3^+$, and suppose that a suffix of ${\mathbf w}$ has the form $h({\mathbf u})$, some ${\mathbf u}\in\Sigma_3^\omega$, where $h=[a,b,c]$. Suppose further that
\begin{enumerate}
\item The word $b$ is a  prefix of $a$, which is a  prefix of $c$.
\item We have $|b|\ge |c|/2$. 
\item The words $b$ and $c$ have a common suffix $s$ such that $|s|\ge |b|/2$.
\end{enumerate}
Then ${\mathbf u}$ doesn't contain any of the factors $00$, $11$, $21$, and $22$.
\end{lemma}
\begin{proof} We treat the factors one at a time. For each factor we rule out 1-letter right extensions of that factor.\vspace{.1in}

\noindent Factor $11$: If $11$ is a factor of ${\mathbf u}$, then so is a factor $v$,  where $v=110$, $v=111$, or $v=112$. In each case, $h(v)$ has the prefix $bbb$ which is a $3$ power. Since ${\mathbf w}$ doesn't contain $3$ powers, the factor $11$ doesn't occur in ${\mathbf u}$.\vspace{.1in}

\noindent Factor $00$: If $00$ is a factor of ${\mathbf u}$, then so is a factor $v$, where $v=000$, $v=001$, or $v=002$. In each case
$h(v)$
has the prefix $aab$ since $b$ is a prefix of $a$ and $c$. However, $aab$ is an $r$ power with $r\ge 5/2$, since $b$ is a prefix of $a$ and $|b|\ge|c|/2\ge|a|/2$. \vspace{.1in}

\noindent Factor $21$: If $21$ is a factor of ${\mathbf u}$, then so is a factor $v$, where $v=211$, $v=210$, or $v=212$. In each case, $h(v)$ has the prefix $cbb$, since $b$ is a prefix of $a$ and $c$. This has the suffix $sbb$, since $s$ is a suffix of $c$. But $sbb$ is an $r$ power with $r\ge 5/2$.\vspace{.1in}

\noindent Factor $22$: If $22$ is a factor of ${\mathbf u}$, then so is a factor $v$, where $v=220$, $v=221$, or $v=222$. In each case,  $h(v)$ has the prefix $ccb$, which is an $r$ power with $r\ge 5/2$, since $b$ is a prefix of $c$ and $|b|\ge|c|/2$.
\end{proof}
It is convenient to work with $\gamma^2=[012,02,0102]$ rather than simply $\gamma$. If $h=[a,b,c]$, then $h\circ\gamma^2=[abc,ac,abac]$, and each image of a letter under $h\circ\gamma^2$ ends in $c$. This allows us to define
$h'=c(h\circ\gamma^2)c^{-1}=[cab,ca,caba]$.
\begin{lemma}\label{iterate}
Suppose that ${\mathbf w}\in\Sigma_3^\omega$ is low.  Suppose that ${\mathbf w}$ has the form $h({\mathbf u})$, some ${\mathbf u}\in\Sigma_3^\omega$ where $h=[a,b,c]$, some $a,b,c\in\Sigma_3^+$ such that
\begin{enumerate}
\item[i.] The word $b$ is a  prefix of $a$, which is a  prefix of $c$.
\item[ii.] We have $|b|\ge |c|/2$. 
\item[iii.] The words $b$ and $c$ have a common suffix $s$ such that $|s|\ge |b|/2$.
\end{enumerate}
Then  a suffix of ${\mathbf w}$ has the form $h(\gamma^2({\mathbf v}))$,  
some ${\mathbf v}\in\Sigma_3^\omega$; equivalently,  
a suffix of ${\mathbf w}$ has the form $h'({\mathbf v}))$ 
some ${\mathbf v}\in\Sigma_3^\omega$, where $h'=c(h\circ\gamma^2)c^{-1}=[cab,ca,caba]$. Letting
$A=h'(0)$, $B=h'(1)$, and $C=h'(2)$, we have
\begin{enumerate}
\item\label{prefix} The word $B$ is a  prefix of $A$, which is a  prefix of $C$.
\item \label{c<2b} We have $|B|\ge |C|/2$. 
\item \label{suffix} The words $B$ and $C$ have a common suffix $S$ such that $|S|\ge |B|/2$.
\end{enumerate}
\end{lemma}
\begin{proof}
Replacing ${\mathbf w}$ by one of its suffixes if necessary, suppose that ${\mathbf w}$ has the form $h({\mathbf u})$.  By Lemma~\ref{00,11,21,22}, words $00$ and $11$ are not factors  of ${\mathbf u}$. Therefore, if $2$ is not a factor of some suffix of ${\mathbf u}$, then that suffix is $(01)^\omega$ or $(10)^\omega$. This is impossible since ${\mathbf w}$ is low. 
We conclude that every suffix of ${\mathbf u}$ contains a $2$.

Parse a suffix of ${\mathbf u}$ into blocks starting with $2$, and containing a single $2$. Words $21$ and $22$ are not factors of ${\mathbf u}$; thus the blocks begin with $20$. Since $00$ and $11$ are not factors, $0$'s and $1$'s alternate in the blocks.  The block $20101$ is impossible; it would imply a block $z=2010101$ or $z=201012$. In both cases, $h(z)$ begins with
$cababa$, since $a$ is a prefix of $c$. However, $ababa$ is an $r$ power with $r\ge 5/2$. It follows that a suffix of ${\mathbf u}$ is concatenated from blocks $20$, $201$, and $2010$. Therefore, a suffix of ${\mathbf w}$ is concatenated from blocks $ca=A$, $cab=B$, and $caba=C$, and has the form $h'({\mathbf v})$ for some ${\mathbf v}\in\Sigma_3^\omega$.

We see that
property \ref{prefix} certainly holds. Also, 
$$|C|=|caba|\le 2|ca|=2|B|$$
since $|c|\ge|b|$, establishing \ref{c<2b}.

Finally, the words $C$ and $B$ have the common suffix $S=sa$. However,
$$2|s|+|a|\ge|b|+|a|\ge 2|b|\ge |c|.$$
Thus
$$|S|=|sa|\ge|c|-|s|=|ca|-|sa|=|B|-|S|,$$
so that $2|S|\ge|B|$, establishing \ref{suffix}. 
\end{proof}

\begin{proof}[Proof of morphism lemma] 
Define morphisms $h_n$ recursively by $h_1=h$, and $h_{n+1} =ch_n\circ\gamma^2c^{-1}$. By induction on the previous lemma, for each positive integer $n$, a suffix of ${\mathbf w}$ has the form $h_n({\mathbf v})$,  
some ${\mathbf v}\in\Sigma_3^\omega$. 

Having a suffix of the form $h_{n+1}({\mathbf v})$ is equivalent to having a suffix of the form $h_n(\gamma^2
({\mathbf v}))$. By induction on $n$, having a suffix of
the form $h_n({\mathbf v})$ is equivalent to having a suffix of the form $h(\gamma^{2n}({\mathbf v}))$. Thus a suffix of ${\mathbf w}$ has the form $h(\gamma^{2n}({\mathbf v}))$. It follows that whenever $g$ is a prefix of ${\mathbf G}$, then $h(g)$ is a factor of ${\mathbf w}$, and the lemma follows.
\end{proof}

\begin{proof}[Proof of dual morphism lemma]
Replacing ${\mathbf z}$ by a suffix if necessary, suppose that ${\mathbf z}$ has the form $h({\mathbf u})$.
Let $T=\{z\in\Sigma_3^*:h(z)\text{ is a factor of }{\mathbf z}\}.$ 
Let $T^R=\{z^R:z\in S\}$. Then $T^R$ is an right factorial language. Using K\"{o}nig's lemma, choose a word ${\mathbf z}'\in\Sigma_3^\omega$ such that every prefix of ${\mathbf z}'$ is in $T^R$. 

Define $h^R=[h(0)^R,h(1)^R,h(2)^R]$.
Let ${\mathbf Z}=h^R({\mathbf z}').$
Then ${\mathbf w}={\mathbf Z}$ and $h'=h^R$ satisfy the conditions of the morphism lemma, so that for every
$g\in{\mathbf G}$, word $h^R(g)$ is a factor of ${\mathbf Z}
=h^R({\mathbf z}')$.

Let $p$ be a prefix of ${\mathbf z}'$ such that 
$h^R(g)$ is a factor of $h^R(p).$ Since every prefix of 
${\mathbf z}'$ is in $T^R$, write $p=z^R$ where $z\in T$.
Then $h^R(g)$ is a factor of $h^R(z^R)$. Thus $h(g^R)$ is a factor of $h(z)$.
By the definition of $T$, 
$h(z)$ is a factor of ${\mathbf z}$, so that $h(g^R)$ is a factor of ${\mathbf z}$.
\end{proof}
\begin{proof}[Proof of first parsing lemma] Let ${\mathbf w}$ be an eligible right-infinite word which omits the words of $S$. Every factor of ${\mathbf w}$ of length $d+1$ is low, and omits words of $S$. By the condition of the lemma it cannot omit $b$. Thus $b$ appears infinitely often in ${\mathbf  w}$, so that a suffix can be concatenated from the $b$-blocks of ${\mathbf w}$. These blocks   
have the form $bu$ where
\begin{itemize}
\item The word $bub$ is a factor of ${\mathbf w}$;
\item Whenever $p$ is a prefix of $bub$ such that $b$ is a suffix of 
$p$, then $p=b$ or $p=bub$.
\end{itemize}
Such a word $u$ is a factor of ${\mathbf w}$ and is therefore low and omits words of $S$. If $|u|>d$ it must therefore have $b$ as a factor. In this case, write $u=vbz$. Then $p=bvb$ is a prefix of $bub$ distinct from $b$ and  $bub$. This is a contradiction.
\end{proof}
\begin{proof}[Proof of second parsing lemma]
Write ${\mathbf w}=b_1b_2b_3b_4\cdots$ such that the $b_i\in B$. Then for $j\ge 2$, if $b_j=s$ let $b_{j-1}=r$ and $b_{j+1}=t$. Since  $rst$ is a factor of ${\mathbf w}$ it is low and cannot contain a factor from $S$. This is a contradiction. We conclude that the suffix $b_2b_3b_4\cdots$ of ${\mathbf w}$ is concatenated from word of $B-\{s\}$.
\end{proof}
\section{A sharpened structure theorem and open problems}
Consider the morphism which is the reverse of $\gamma$, namely
$\gamma^R=[10,2,20]$. Let ${\mathbf G}'$ be its fixed point $${\mathbf G}'=(\gamma^R)^\infty(2)=20102102010\cdots.$$
(An anonymous referee points out that ${\mathbf G}'=\pi({\mathbf p})$ where ${\mathbf p}$ is the sequence studied in  \cite{currie2022} and $\pi$ is the permutation of $[2,0,1]$ of $\Sigma_3$. See the entry for
sequence \href{https://oeis.org/A287072}{A287072} in the On-Line Encyclopedia of Integer Sequences.)
Then the factors of ${\mathbf G}'$ are precisely the reverses of the 
factors of ${\mathbf G}$. One can easily adjust the proof of the 
morphism lemma to prove
\begin{lemma}[Alternate Dual Morphism Lemma]\label{alternate dual morphism lemma}
Suppose that ${\mathbf z}\in\Sigma_3^\omega$ is low.  Suppose that a suffix of ${\mathbf z}$ has the form $h({\mathbf u})$, some ${\mathbf u}\in\Sigma_3^\omega$ where $h=[a,b,c]$, some $a,b,c\in\Sigma_3^+$ such that
\begin{enumerate}
\item The word $b$ is a  suffix of $a$, which is a  suffix of $c$.
\item We have $|b|\ge |c|/2$. 
\item The words $b$ and $c$ have a common prefix $p$ such that $|p|\ge |b|/2$.
\end{enumerate}
Then if $g$ is any factor of ${\mathbf G}'$, then $h(g)$ is a factor of ${\mathbf z}$.
\end{lemma}

As remarked earlier, it is possible to keep track of the morphisms used in the proof of the structure theorem. Doing so gives the following:

\begin{theorem}[Sharpened Structure Theorem]\label{sharpened structure theorem} Let ${\mathbf w}\in\Sigma_3^\omega$ be eligible. One of the following holds.
\begin{itemize}
\item 
For some non-erasing morphism $h$ in Figure~\ref{morphisms} and a permutation $\pi$ of $\Sigma_3$, the word $\pi(h(g))$ is a factor of ${\mathbf w}$ whenever $g$ is a factor of ${\mathbf G}$;
\item 
For some non-erasing morphism $h$ in Figure~\ref{dual morphisms} and some permutation $\pi$ of $\Sigma_3$, the word $\pi(h(g))$ is a factor of ${\mathbf w}$ whenever $g$ is a factor of ${\mathbf G}'$.
\end{itemize}
\end{theorem}

\begin{figure}
\begin{center}
\begin{tabular}{|l|l|l|}\hline
$h(0)$&$h(1)$&$h(2)$\\\hline
$1210202$& $12102$& $1210202102$\\\hline
$1210202102$& $1210202$& $121020210202$\\\hline
$021201212$& $0212012$& $021201212012$\\\hline
$21021102211021102$& $2102110221102$& $2102110221102110221102$\\\hline
\end{tabular}
\end{center}
\caption{Potential morphisms $h$ such that some eligible word contains $h(g)$ for each factor $g$ of ${\mathbf G}$.}\label{morphisms}
\end{figure}

\begin{figure}
\begin{center}
\begin{tabular}{|l|l|l|}\hline
$h(0)$&$h(1)$&$h(2)$\\\hline
$10212$& $0212$& $0210212$\\\hline
$2120210$& $20210$& $202120210$\\\hline
$121201210$& $1201210$& $120121201210$\\\hline
$1201212010$& $1212010$& $121201212010$\\\hline
$1212010$& $12010$& $1201212010$\\\hline
$221102210$& $21102210$& $2110221102210$\\\hline
\end{tabular}
\end{center}
\caption{Potential morphisms $h$ such that some eligible word contains $h(g)$ for each factor $g$ of ${\mathbf G}'$.}
\label{dual morphisms}
\end{figure}

The morphisms in Figures \ref{morphisms} and \ref{dual morphisms} were output by the computer implementation of our proof. In fact, the implementation gave 26 different morphisms, each produced by applying a permutation to one of the listed morphisms. Some further simplification is possible. For example, consider the third morphism in Figure 1, namely
$$h=[021201212,0212012,021201212012].$$ This is a conjugate of
$$h_1=[201212021,2012021,201212012021].$$ Applying the permutation $\sigma=[1,2,0]$ we find
$$\sigma(h_1)=[012020102,0120102,012020120102]=\gamma^4.$$
Morphism $h$ could thus be replaced in Figure~\ref{morphisms} by the identity
morphism. We did not complicate our proof implementation by searching for simplifications/duplications of morphisms.
\begin{problem}Give a version of the sharpened structure theorem with an `aesthetically optimal' list of morphisms.
\end{problem}

The sharpened theorem says that the set of factors of ${\mathbf w}$  
is the same as the set of factors of $h({\mathbf G})$ for a morphism in Figure~\ref{morphisms}, or $h({\mathbf G}')$ for a morphism in Figure~\ref{dual morphisms}. It is not known for which, if any, of the $h$  the words $h({\mathbf G})$ and $h({\mathbf G}')$ have factor complexity $2n+1$, and what the critical exponents of these words are. We make the following conjecture:

\begin{conjecture} For each morphism $h$ in Figure~\ref{morphisms}, the word $h({\mathbf G})$ has factor complexity $2n+1$ and critical exponent $\mu$. For each morphism $h$ in Figure~\ref{dual morphisms}, the word $h({\mathbf G}')$ has factor complexity $2n+1$ and critical exponent $\mu$.
\end{conjecture}

Taken together with 
Carpi and De Luca's result \cite{luca} for Sturmian words, we have the beginnings of a Dejean's Theorem for Arnoux-Rauzy words; we know the minimal critical exponent in the cases $k=2$ and $k=3$.

\begin{problem}For each integer $k\ge 4$, give an Arnoux-Rauzy word over $\Sigma_k$ with the least critical exponent.
\end{problem}
\section*{Appendix: `Human readable' proof of the resolution of the case $S=\{01,12,20,00\}$ of the structure theorem}
The Python implementation outputs the following proof:\vspace{.1in}\\
We assume $F$ includes  [01, 12, 20, 00].\\

Every good word longer than  25 must include factor(s)\\

[0221, 1021, 1022, 1102, 2102, 2110].\vspace{.1in}\\

{\em Author's comment: In step B of the algorithm, we found a length-250 low word $w$ omitting $\{01,12,20,00\}$. In step C of the algorithm, for $n=4$, we find that $w$ has more than $2n+1$ length-4 factors. Of these, 6 (announced above) are found to be needed factors in step C.b of the algorithm.}\vspace{.1in}\\ 

The word w must omit a factor from\\

[0210, 0211, 2210, 2211].\vspace{.1in}\\

{\em Author's comment: These are the 4 unneeded factors found in step C.b of the algorithm. They cannot all be included in ${\mathbf w}$, because $6+4=10$ would be too many length-4 factors in ${\mathbf w}$.}\vspace{.1in}\\
This gives rise to  4  cases.\vspace{.05in}\\
    Case 1: $F$ includes  [01, 12, 20, 00, 0210]\\
    Case 2: $F$ includes  [01, 12, 20, 00, 0211]\\
    Case 3: $F$ includes  [01, 12, 20, 00, 2210]\\
    Case 4: $F$ includes  [01, 12, 20, 00, 2211]\vspace{.1in}\\

{\em Author's comment: This is the case division of step C.b.ii.}\vspace{.1in}\\

Case 1: $F$ includes  [01, 12, 20, 00, 0210].\\

Every good word longer than  68 must include factor(s)\\

[02110, 02210, 02211, 10211, 10221, 11022, 21021, 21102, 22102, 22110].

The word w must omit a factor from\\

[11021, 21022].

This gives rise to  2  cases.\vspace{.05in}\\
Case 1.1: $F$ includes  [01, 12, 20, 00, 0210, 11021]\\
Case 1.2: $F$ includes  [01, 12, 20, 00, 0210, 21022]\\

Case 1.1: $F$ includes  [01, 12, 20, 00, 0210, 11021].\\

Every good word longer than  81 must include factor(s)\\

[02110221022, 02110221102, 02210211022, 02210221102, 02211022102, 10211022102, 10211022110, 10221021102, 10221022110, 10221102210, 11022102110, 11022102211, 11022110221, 21021102210, 21021102211, 21022110221, 21102210211, 21102210221, 21102211022, 22102110221, 22102211022, 22110221021, 22110221022].

In this case w must omit the factor  02110221021.\vspace{.1in}\\

{\em Author's comment: Here, because only 1 unneeded factor was found  in Step C.b, we don't make a new subcase.}\vspace{.1in}\\
    The $210$-blocks  are among  [2102, 21021102, 210221102, 2102110221102].\vspace{.1in}\\

{\em Author's comment: It turns out that when we find a length-250 low word $w$ omitting the current list of missing factors, the complexity is $2n+1$ for $n=1,2,\ldots, 20$. We therefore, in step C.a, search for a morphism, considering factors $b$ of $w$ with $|b|\le 3$ and finding the $b$-blocks of $w$.}\vspace{.1in}\\

    The word  2102 cannot occur twice in a concatenation of these\vspace{.1in}\\

{\em Author's comment: If 2102 occurs twice, it occurs in the context $s2102t$, where $s,t$ are two of the $210$-blocks. However, one checks that each such word $s2102t$ contains an $r$ power with $r\ge 5/2$.}\vspace{.1in}\\
    Thus a suffix of w is  concatenated from  [21021102, 210221102, 2102110221102].\\

    Taking conjugates, a suffix is the image under morphism\\

[221102210, 21102210, 2110221102210].

    which satisfies the dual morphism lemma.\vspace{.1in}\\

(This concludes Case 1.1).\\
{\em Author's comment: Since ${\mathbf w}$ has a suffix of the form $h({\mathbf u})$, the dual morphism lemma shows that
${\mathbf w}$ has critical exponent at least $\mu$, and the current subcase is resolved.}\vspace{.1in}\\

Case 1.2: $F$ includes  [01, 12, 20, 00, 0210, 21022].\\

    The 210-blocks  are\\

[21021102, 2102110221102, 21021102211021102, 2102110221102110221102].

    The word  21021102 cannot occur twice in a concatenation of these.\\

    Thus a suffix of w is concatenated from\\

[2102110221102, 21021102211021102, 2102110221102110221102]

    Taking conjugates, a suffix is the image under morphism\\

[21021102211021102, 2102110221102, 2102110221102110221102].

    which satisfies the morphism lemma.\\

(This concludes Case 1.2).\\

(This concludes Case 1).\\

Case 2: $F$ includes  [01, 12, 20, 00, 0211].\\

Every good word longer than  54 must include factor(s)\\

[02102, 02210, 02211, 10210, 10221, 11021, 11022, 21021, 21022, 21102, 22102, 22110].

The word w has  $12 > 2 \times 5 + 1$ factors of length  5. This is impossible.\\

(This concludes Case 2).\\
\vspace{.1in}\\
{\em Author's comment: We found too many needed factors in step 
C.b.}\vspace{.1in}\\

Case 3: $F$ includes  [01, 12, 20, 00, 2210].\\

Every good word longer than  56 must include factor(s)\\

[02102, 02110, 02211, 10210, 10211, 10221, 11021, 11022, 21021, 21022, 21102, 22110].

The word w has  $12 > 2 \times 5 + 1$ factors of length  5. This is impossible.\\

(This concludes Case 3).\\

Case 4: $F$ includes  [01, 12, 20, 00, 2211].\\

Every good word longer than  41 must include factor(s)\\

[02102, 02110, 02210, 10210, 10211, 10221, 11021, 11022, 21021, 21022, 21102, 22102].

The word w has  $12 > 2 \times 5 + 1$ factors of length  5. This is impossible.\\

(This concludes Case 4).\vspace{.1in}\\
{\em Author's comment: The proof in the present case is much shorter than for the resolution of $S=\{11,22\}$. In that longer resolution the computer frequently invokes step A of the algorithm (not used in the present case) to avoid repeating case analysis which was already carried out.}


\begin{thebibliography}{24}
\bibitem{allouche}
J.-P. Allouche and J. Shallit, {\em Automatic sequences: theory, applications, generalizations}, Cambridge, Cambridge University Press, 2003.
\bibitem{andrade}
Jonathan~Andrade, 
Avoiding additive powers in words, 
undergraduate Honour's thesis, Thompson Rivers University, 2024. Available at \url{https://tru.arcabc.ca/islandora/object/tru%3A6415}.

 \bibitem{baranwal}Aseem Baranwal, James Currie, Lucas Mol, Pascal Ochem, Narad Rampersad, and Jeffrey O. Shallit, Antisquares and critical exponents, {\em Discrete Math.
 \& Theoret. Comput. Sci.} {\bf 25:2} \#11 (2023).

\bibitem{carpi}
Arturo Carpi,
On Dejean’s conjecture over large alphabets,
{\em Theoret. Comp. Sci.} {\bf 385}, 137--151 (2007).

\bibitem{luca}Arturo Carpi and Alessandro De Luca, Special factors, periodicity, and an application to Sturmian words. {\em Acta Inform.} {\bf 36}, 983--1006 (2000).

\bibitem{cassaigne} Julien~Cassaigne, James~D.~Currie, Luke~Schaeffer, and Jeffrey~Shallit,
Avoiding three consecutive blocks of the same size and same sum,
{\em J. ACM} {\bf 61}(2), 1--17 (2014).

\bibitem{cassaigne2017} Julien Cassaigne, S\'{e}bastien Labb\'{e}, and Julien Leroy, A set of sequences of complexity $2n+1$, \url{https://arxiv.org/abs/1707.02741}.

\bibitem{morteza} James D. Currie and M. Mohammad-Noori,
Dejean’s conjecture and Sturmian words
{\em Eur. J. Combin.} {\bf 28}, 876--890 (2007).
 
\bibitem{mol1} James D. Currie, Lucas Mol, and Narad Rampersad, The repetition threshold for binary rich words. {\em Discrete Math.
 \& Theoret. Comput. Sci.} {\bf 22}(1), DMTCS–22–1–6 (2020). 

\bibitem{mol2} James D. Currie, Lucas Mol, and Jarkko Peltom\"aki
, The repetition threshold for ternary rich words. {\em Electron. J. Comb.} {\bf 32}(2), P2.55 (2025). 

\bibitem{currie2022} James Currie, Pascal Ochem, Narad Rampersad and Jeffrey Shallit, Properties of a ternary infinite word, {\em RAIRO-Th\'{e}or. Inf. Appl.} {\bf 57} (2023). 

\bibitem{dekking}
James~Currie, Lucas~Mol, Narad~Rampersad, and Jeffrey~Shallit,
Extending Dekking's construction of an infinite binary word avoiding abelian 4-powers,
{\em SIAM J. Discrete Math.} {\bf 38}, 2913--2925 (2024).

\bibitem{template}
James~D.~Currie and Narad~Rampersad,
Fixed points avoiding Abelian $k$-powers,
{\em J. Comb. Theory Ser. A} {\bf 119}(5), 942--948 (2012).

\bibitem{rampersad1}
James Currie and Narad Rampersad,
A proof of Dejean's conjecture,
{\em Math. Comp.} {\bf  80}, 1063--1070 (2011).

\bibitem{dejean}
Fran\c{c}oise Dejean,
Sur un th\'eor\`eme de Thue,
{\em J. Combin. Theory Ser. A} {\bf  13}, 90--99 (1972).

\bibitem{dvorakova2} L\'ubomíra ~Dvo\v{r}\'akov\'a, Kate\v{r}ina~Medkov\'a, and Edita~Pelantov\'a, 
Complementary symmetric Rote sequences: the critical exponent and the recurrence function,
 {\em Discrete Math. \& Theoret. Comput. Sci.}  DMTCS--22--1--20  (2020).

\bibitem{dvorakova1}
L\'ubomíra ~Dvo\v{r}\'akov\'a, Daniela~Opo\v{c}ensk\'a, Edita~Pelantov\'a, and Arseny~M. Shur,
 On minimal critical exponent of balanced sequences, {\em Theoret. Comput. Sci.}
{\bf 922}, 158--169 (2022).

\bibitem{kleene} Stephen Kleene, {\em  Mathematical Logic}, New York, Wiley, 1967.

\bibitem{rosenfeld1}
Florian~Lietard and Matthieu~Rosenfeld,
Avoidability of additive cubes over alphabets of four numbers,  
in N.~Jonoska and D.~Savchuk, editors, {\em Developments in Language Theory 2020}, Vol.~12086 of {\em Lecture Notes in Computer Science}, pp.~192--206, Springer-Verlag, 2020.

\bibitem{mousavi} Hamoon Mousavi, Automatic theorem proving in Walnut. Preprint: \url{https://arxiv.org/abs/1603.06017}(2016, updated 2021).

\bibitem{ollagnier}
Jean Moulin~Ollagnier,
Proof of Dejean’s conjecture for alphabets with 5,6,7,8,9,10
 and 11 letters, {\em Theoret. Comp. Sci.} {\bf 95}, 187--205 (1992).

\bibitem{pansiot}
Jean-Jacques Pansiot, \`A propos d’une conjecture de F. Dejean sur les r\'ep\'etitions dans les mots,
{\em Disc. App. Math.} {\bf 7}, 297--311 (1984). 

\bibitem{rampersad2} Narad Rampersad, Jeffrey Shallit, and \'Elise Vandomme, Critical exponents of infinite balanced words. {\em Theoret. Comput. Sci.} {\bf 777}, 454--463 (2020).

\bibitem{rao}
Mich\"ael Rao,
Last cases of Dejean's conjecture,
{\em Theoret. Comput. Sci.} {\bf 412}, 3010--3018 (2011).

\bibitem{rosenfeld2}
Mich\"ael~Rao and Matthieu~Rosenfeld,
Avoiding two consecutive blocks of same size and
same sum over $\mathbb{Z}^2$, 
{\em SIAM J. Discrete Math.}
{\bf 32}(4) (2018), 2381--2397.

\bibitem{shallit2022}
Jeffrey~Shallit,
{\em The Logical Approach To Automatic Sequences: Exploring
  Combinatorics on Words with {\tt Walnut}}, Vol.~482 of {\em London Math. Soc.
  Lecture Note Series}, Cambridge University Press, 2022.

\bibitem{walnut}Jeffrey O. Shallit, Walnut papers and books: \url{https://cs.uwaterloo.ca/~shallit/walnut-papers.html} (2025).

\bibitem{shallit2019}Jeffrey Shallit and Arseny Shur, Subword complexity and power avoidance, {\em Theoret. Comput. Sci.} {\bf 792}, 96--116  (2019).

\end{thebibliography}
\end{document}